\documentclass{article}

\usepackage{amsmath}
\numberwithin{equation}{section}
\usepackage{amssymb}
\usepackage{theorem}
\usepackage{enumitem}

\usepackage[english]{babel}
\usepackage[utf8]{inputenc}

\theoremstyle{plain} {\theorembodyfont{\slshape}
\newtheorem{thm}{Theorem}[section]

}

{\theorembodyfont{\rmfamily}
\newtheorem{exa}[thm]{Example}

\newtheorem{rem}[thm]{Remark}
}

\newtheorem{lem}[thm]{Lemma}

\newcommand{\proof}{{\bf Proof:\ }}
\newcommand{\Endproof}{\hspace*{\fill} $\Box$ \vspace{1ex} \noindent }

\newcommand{\X}{\textbf{X}}

\newcommand{\Z}{\textbf{Z}}
\newcommand{\Sum}{\textbf{S}}
\newcommand{\thetav}{\boldsymbol\theta}
\newcommand{\x}{\textbf{x}}
\newcommand{\y}{\textbf{y}}
\newcommand{\p}{\textbf{p}}
\newcommand{\cb}{\textbf{c}}
\newcommand{\ab}{\textbf{a}}
\newcommand{\nolla}{\textbf{0}}

\newcommand{\RR}{\mathbb R}
\newcommand{\Ss}{\mathbb S}

\newcommand{\hanta}{\overline{F}}
\newcommand{\PP}{\mathbb P}
\newcommand{\EE}{\mathbb E}

\newcommand{\dist}{\textup{dist}}
\newcommand{\ud}{\mathrm{d}}

\usepackage[T1]{fontenc} 
\usepackage{setspace}

\title{Precise asymptotics of ruin probabilities for a class of multivariate heavy-tailed distributions}

\author{Miriam Hagele}

\begin{document}
\maketitle

\begin{abstract}
\noindent
This article studies asymptotic approximations of ruin probabilities of multivariate random walks with heavy-tailed increments. Under our assumptions, the distributions of the increments are closely connected to multivariate subexponentiality and admit dependence between components.
\end{abstract}

\section{Introduction}
Insurance companies operating in different regions or offering different types of insurances can use multivariate models to estimate their risks and to find suitable hedging strategies.
In this approach, lines of business are considered as components of multivariate random vectors. 
Often, these models are built on the assumptions that the risks follow multivariate regularly varying distributions. One example is in \cite{Resnick1}. 
However, the risk process of an insurance company might not always be regularly varying but instead be characterized by a subexponential distribution, see \cite{Asmussen2} or \cite{Baltrunas1}. On the real line, the distribution $F$ of a positive random variable is subexponential if $\lim_{x\to\infty} \overline{F^{*2}}(x)/\overline{F}(x) =2$ where $F^{*2}$ denotes the convolution of $F$ with itself. We are interested in multivariate generalizations of this concept in the setting of random walks.

Mathematically, we examine the multivariate random walk
\begin{displaymath}
\Sum_0=\nolla,\qquad \Sum_n=\X_1+\dots+\X_n, \quad n\geq 1,
\end{displaymath}
where $\X,\X_1,\X_2,\dots\in \RR^d$ are i.i.d. random vectors on a probability space $(\Omega,\mathcal{F},\PP)$. To simplify notations, we use polar coordinates. 
We write $\X=R\thetav$, where $R= \|\X\|$ is a subexponentially distributed real-valued random variable describing the length of the vector $\X$ and $\thetav=\X/\|\X\|\in \Ss^{d-1}$ is a random vector on the unit sphere representing the angle or direction of $\X$. 
We denote by $\|\cdot\|$ the $L_1$-norm and the set $\Ss^{d-1}$ denotes the $L_1$ unit sphere or diamond. The set $\Ss^{d-1}_{>0}:= \{\x\in\Ss^{d-1}:x^1>0,\dots,x^d>0\}$ is the subset of $\Ss^{d-1}$ where all components are positive.
The random vector $\X$ can be interpreted as the yearly net-payout of an insurance company. Due to this interpretation, we assume below that all components of $\X$ have a negative expectation i.e.\ all lines of business are on average profitable.

The objective of this article is to study the probability that the random walk $\{\Sum_n\}$ hits, at some time $n$, the set $uC_\delta$ defined as 
\begin{equation} \label{eq0cdelta}
uC_{\Theta,\delta}= uC_\delta =\left\{\x\in \RR^d:\|\x\|>u,\frac{\x}{\|\x\|}\in \Theta^\delta\right\},
\end{equation}
where $\Theta^\delta$ is the $\delta$-swelling of the set $\Theta$, specified in Assumption \ref{as0thetav}, in the set $\Ss^{d-1}_{>0}$. We define the $\delta$-swelling in $\Ss^{d-1}_{>0}$ by 
\begin{displaymath}
\Theta^\delta :=\{\x\in \Ss^{d-1}_{>0}: \|\x-\y\|<\delta \textrm{ for some }\y \in \Theta\},
\end{displaymath}
see \cite{Das2}.
Hence, the set $uC_\delta$ is a truncated cone in the positive orthant. It consists of all vectors in the direction of the set $\Theta^\delta$ with $L_1$-length greater than some threshold $u$.

Ruin probabilities for similar sets in the case where $\X$ follows a multivariate regular varying distribution are derived in \cite{Hult1}. Different concepts of multivariate subexponentiality were introduced in \cite{Cline1}, \cite{Omey1} and \cite{Samorodnitsky}. We will use the latter approach. \cite{Samorodnitsky}, the closest reference to this article, studies ruin probabilities in sets where at least one component has to be large. Here, we show that in the case of asymptotically dependent components the ruin probability of sets that contain $\Theta$, the asymptotic support of the random vector $\theta$, is asymptotically equivalent to the integrated tail distribution of the vector length.

Throughout the article, bold letters denote vectors, upper case and Greek letters are for random vectors and lower case for non-random vectors. The components are denoted by upper indices, for instance, $\x=(x^1,\dots,x^d)$, $\nolla=(0,\dots,0)$. Relation symbols, such as $\x<\y$ are meant component-by-component. The notation $f\in o(g)$ means $\lim_{u\to\infty} (f(u)/g(u)) = 0$ and $f\sim g$ refers to $f(u)/g(u) \to 1$ as $u\to\infty$. Unless otherwise specified, the asymptotic relation holds for $u\to\infty$.
The tail distribution of a distribution function $F$ is denoted by $\hanta = 1-F$ and the integrated tail function by $\hanta_I(u)=\min\left(1,\int_u^\infty \hanta(v)\ud v\right)$. The set $B(\x,r)$ defines a ball in $L_1$-norm with radius $r$ centred at the point $\x$ and $\textrm{cl}(A)$ denotes the closure of a set $A$.

\section{Main results}
The main results illustrate the asymptotic behaviour of a random walk consisting of the sum of i.i.d random vectors $\X=R\thetav$, where $R$ represents the length of the random vector and $\thetav$ indicates its direction. Setting $\cb:=-\EE(\X)$ and
\begin{equation} \label{eq0seta}
A:=\left\{\x\in\RR^d:\sum_{k=1}^d x^k >1\right\},
\end{equation} 
the assumptions on $\X$ are:
\begin{enumerate}[label=(A\arabic*)]
\item \label{as0r}
$R$ is subexponentially distributed with distribution function $F$.
\item \label{as0inttail}
$\hanta_I(\gamma u)=o(\hanta_I(u))$ for all $\gamma>1$.
\item \label{as0thetav}
There exists a set $\Theta$ with $\textrm{cl}(\Theta) \subset \Ss^{d-1}_{>0}$ such that 
\begin{equation} \label{eq0thetav}
\lim_{h\to\infty}\PP\left(\thetav \in \Theta^\varepsilon  |R>h\right) = 1  \quad\textrm{for all }\varepsilon>0.
\end{equation}
\item \label{as0expectation}
$-\EE(\X) >\nolla$.
\item \label{as0finite}
$\int_0^\infty \PP(\X\in A+v\cb)\ud v<\infty$.
\item \label{as0subexp}
The distribution defined by $H(u)=\max\left(0,1-\int_0^\infty \PP(\X\in uA+v\cb)\ud v\right)$ is subexponential.
\end{enumerate}

The assumptions \ref{as0r} and \ref{as0inttail} imply that $R$ is a heavy-tailed distribution with a lighter tail than regularly varying distributions. For instance, Weibull and lognormal distributions that possess all their power moments satisfy these assumptions. 
Assumption \ref{as0thetav} excludes the possibility of asymptotically independent components. One advantage of the assumption is the possibility to predict upper and lower bounds for the magnitude of losses in the other components if one knows a large outcome in one component. These kind of situations arise for instance if the reporting time of events varies between different lines of business. After the first line of business reports a large observation, one can compute in which range the magnitude of losses are in the other lines of business. In the bivariate case, if $\Theta$ consist only of one point in $\Ss^1_{>0}$ the components exhibit full asymptotic dependence, see \cite{Resnick2}. If $\Theta$ is an interval in $\Ss^1_{>0}$ the components exhibit strong asymptotic dependence. These concepts come from the case of multivariate regular variation where one looks at the limit measure of the random variable as in \cite{Das1}. Here, we define the dependence structure through the set $\Theta$ for which the conditional probability stated in Assumption \ref{as0thetav} converges to one. 
Assumption \ref{as0expectation} corresponds to a positive safety loading. Similarly to the one-dimensional case where $1-H(u)$ is the integrated tail function, one can interpret $1-H(u)$ in Assumption \ref{as0subexp} as an integrated tail function integrating along the set $uA+v\cb$.

\begin{thm} \label{lem0thetavakio}
Set $\delta>0$ such that for all $\y\in \Theta^\delta$ for all components $y^k> \delta /(4+\delta)$. Set $uC_\delta$ as in (\ref{eq0cdelta}) and assume \ref{as0r} - \ref{as0subexp} hold. Then,
\begin{equation} \label{eq0ruinprob}
\lim_{u\to\infty}\frac{\PP(\Sum_n\in uC_\delta \textrm{ for some }n\geq 1)}{\frac{1}{\|\cb\|}\int_u^\infty \PP(R>v)\ud v}=1.
\end{equation}
\end{thm}

The following theorem states a partially reverse result of the theorem above.

\begin{thm} \label{lem0equivalence}
Assume \ref{as0r}, \ref{as0inttail}, \ref{as0expectation} - \ref{as0subexp} hold, $\lim_{h\to\infty} \PP(\thetav \in \Ss^{d-1}_{>0} |R>h)=1$, the limit  in (\ref{eq0thetav}) exists for all $\Theta^\delta\subset \Ss^{d-1}_{>0}$ and additionally that $u^2\hanta(u)=o(1)$. Then, there exists some set $\Theta$ for which Equation (\ref{eq0ruinprob}) is equivalent to Condition (\ref{eq0thetav}).  
\end{thm}

\begin{rem}
The set $\Theta$ does not have to be the smallest possible set for which Assumption \ref{as0thetav} holds. Choosing the set as small as possible improves the result. 
\end{rem}

\subsection{Preliminary results}
The proof of Theorem \ref{lem0thetavakio} uses similar techniques as the main reference \cite{Samorodnitsky} which applies some ideas of the proof in the one-dimensional counterpart studied in \cite{Zachary}. To show the asymptotic lower bound, we use a novel geometric approach stated in Lemma \ref{lem0geomapulem} while the asymptotic upper bound uses the following remark which is based on Theorem 5.2 of \cite{Samorodnitsky}.

\begin{rem} \label{rem0samsun}
Theorem 5.2 in \cite{Samorodnitsky} shows the asymptotic ruin probability in a model where the claim size vectors have positive increments. In this model, $\Z_i$ are positive $d$-dimensional i.i.d random vectors with the common distribution $G_{\Z}$ and the interarrival times $(Y_i)_{i\geq 1}$ form a sequence of i.i.d positive random variables. Setting $A$ as in (\ref{eq0seta}) and $\X_i=\Z_i-Y_i\p$ where $\p$ is some positive $d$-dimensional vector, we assume $\cb = -\EE(\X) >\nolla$ and that the distribution $\int_0^\infty G_{\Z}(uA+v\cb) \ud v / \int_0^\infty G_{\Z}([0,\infty)^d +v\cb) \ud v$
is subexponential. Then,
\begin{displaymath}
\PP\left(\sum_{i=1}^n \X_i \in uA \textrm{ for some }n\geq 1\right)\sim \int_0^\infty G_{\Z}(uA+v\cb) \ud v \qquad \textrm{ as $u\to\infty$.}
\end{displaymath}
\end{rem}

The geometric approach used in the proof of the asymptotic lower bound is based on the principle of a single big jump. Thus, the main contribution for the ruin probability comes from the case where the random walk jumps from the neighbourhood of its expectation to the ruin set in one step. The following lemma shows that such jump in the direction of the set $\Theta^{\frac{\delta}{2}}$ always hits the set $uC_\delta$ if the length of the jump is large enough.

\begin{lem} \label{lem0geomapulem}
Let $\varepsilon>0$ and $K\in \RR_+$. Assume $\cb>\nolla$ is a vector in $\RR^d$, $\Theta$ a set with $\textrm{cl}(\Theta)\subset\Ss^{d-1}_{>0}$ and choose $\delta>0$ such that $y^k>\frac{\delta}{4+\delta}$ for all $\y\in \Theta^\frac{\delta}{2}$ and $k=1,\dots,d$. Define $B_n:=B(-n\cb,n\varepsilon)$
and 
\begin{displaymath}
u_n := \max\left(u+n\|\cb\|+nd\varepsilon,\frac{4+\delta}{\delta}(n\|\cb\|+nd\varepsilon)\right).
\end{displaymath}
Then for all $n$, $\x+t\y \in uC_\delta$ uniformly for all $\x\in B_n, \y\in \Theta^\frac{\delta}{2}$ and $t>u_n$. The same result holds if we define $B_n:=B(-n\cb,K+n\varepsilon)$ and $u_n:=\max\left(u+n\|\cb\|+nd\varepsilon+dK,\frac{4+\delta}{\delta}(n\|\cb\|+nd\varepsilon+dK)\right)$.
\end{lem}
\proof
Assuming $\x\in B_n$, $\y \in \Theta^\frac{\delta}{2}$ and $t> u_n$, we will show
\begin{displaymath}
\|\x+t\y \|>u \qquad \textrm{ and } \qquad \frac{\x+t\y}{\|\x+ t\y\|}\in \Theta^\delta
\end{displaymath}
to prove $\x+t\y\in uC_\delta$. Taking $\x\in B_n$ implies $x^k\in (-nc^k-n\varepsilon,-nc^k+n\varepsilon)$ for all components and from $y^k>\frac{\delta}{4+\delta}$ for all $k=1,\dots,d$ follows $x^k+ty^k>0$. Therefore, if $t> u+n\|\cb\|+nd\varepsilon$ and $y^k>\frac{\delta}{4+\delta}$ then
\begin{eqnarray*} 
\|\x+t\y\| &=& \sum_{k=1}^d |x^k+ty^k| = \sum_{k=1}^d (x^k+ty^k)\geq \sum_{k=1}^d (-nc^k-n\varepsilon+ty^k)\\
&=& -n\|\cb\|-nd\varepsilon+t\|\y\| > u.
\end{eqnarray*}
Next, we show $\left\|\frac{\x+t\y}{\|\x+t\y\|}-\y\right\|<\frac{\delta}{2}$ to prove that $\frac{\x+t\y}{\|\x+t\y\|} \in \Theta^\delta$. By the calculations above
\begin{eqnarray*} 
\left\|\frac{\x+t\y}{\|\x+t\y\|}-\y\right\| &=&\left\|\frac{\x+t\y-(-\|\x\|+t\|\y\|)\y}{\|\x+t\y\|}\right\|= \frac{\|\x+\|\x\|\y\|}{\|\x+t\y\|} \\
&\leq& \frac{2 \|\x\|}{\|\x+t\y\|}\leq\frac{2n\|\cb\|+2nd\varepsilon}{t\|\y\|-n\|\cb\|-nd\varepsilon}<\frac{\delta}{2}.
\end{eqnarray*}
Similar calculations yield the result for $B_n=B(-n\cb,K+n\varepsilon)$ and its corresponding $u_n$.
\Endproof

\subsection{Proof of Theorem \ref{lem0thetavakio}}
\subsubsection{Lower bound}
Throughout the proofs $\PP(\Sum_n\in uC_\delta \textrm{ for some }n)$ means $n\geq 1$ if not otherwise specified. At first, we show that $\frac{1}{\|\cb\|}\int_u^\infty \PP(R>v) \ud v$ is an asymptotic lower bound of the probability $\PP(\Sum_n\in uC_\delta \textrm{ for some }n\geq 1)$.
By the law of large numbers 
\begin{displaymath}
\PP\left(\Sum_n  \in B(-n\cb,n\varepsilon)\right) \underset{n\to\infty}{\longrightarrow} 1
\end{displaymath}
for all $\varepsilon>0$. Hence, for all $\varepsilon,\delta'>0$ there exists some constants $K:=K_{\varepsilon,\delta'}$ and $n_1$ such that $\PP(\Sum_n\in B_n) > 1-\delta'$ for all $n$ where
\begin{displaymath}
B_n = \left\{\begin{array}{ll} 
B(-n\cb,n\varepsilon+K)& \textrm{if }n\leq n_1\\
B(-n\cb,n\varepsilon)& \textrm{if }n>n_1.\\
\end{array}\right.
\end{displaymath}
Applying Lemma \ref{lem0geomapulem} yields 
\begin{eqnarray} 
&&\PP\left(\Sum_n  \in uC_\delta \textrm{ for some } n\geq 1\right)\nonumber \\
&\geq& \sum_{n\geq 0} \PP\left(\Sum_k  \notin uC_\delta \ \forall k\leq n, \Sum_n \in B_n,\X_{n+1}\in u_n C_{\frac{\delta}{2}}\right)\label{ineq01}\\
&\geq& \sum_{n\geq 0} \left(1-\delta'-\PP\left(\Sum_k \in uC_\delta  \textrm{ for some } k\leq n\right)\right)\PP\left(\X_{n+1}\in u_n C_\frac{\delta}{2}\right)\label{ineq02}\\
&\geq& \left(1-\delta'-\PP\left(\Sum_n \in uC_\delta  \textrm{ for some } n\geq 1\right)\right)\sum_{n\geq 0} \PP\left(\X\in u_n C_\frac{\delta}{2} \right)\nonumber,
\end{eqnarray}
where $u_n$ is defined as in Lemma \ref{lem0geomapulem} depending on the definition of $B_n$.
Inequality (\ref{ineq01}) follows from Lemma \ref{lem0geomapulem} where the random walk makes a big jump from the set $B_n$ to the set $uC_\delta$ in the direction of $\Theta^\frac{\delta}{2}$. Inequality (\ref{ineq02}) is due to the law of large numbers and the fact that $\PP(A\cap B)=\PP(B\backslash A^c)\geq \PP(B)-\PP(A^c)$. Rearranging this inequality yields for large $u$
\begin{eqnarray*} 
\PP\left(\Sum_n \in uC_\delta \textrm{ for some }n\geq 1\right) 
&\geq& \frac{(1-\delta')\sum_{n\geq 0}\PP\left(\X\in u_n C_\frac{\delta}{2}\right)}{1+\sum_{n\geq 0}\PP\left(\X\in u_n C_\frac{\delta}{2}\right)} \\
&\geq& \frac{1-\delta'}{1+\delta''}\sum_{n\geq 0}\PP\left(\X\in u_n C_\frac{\delta}{2}\right)
\end{eqnarray*}
for any given $\delta''>0$.
Since $\PP\left(\X\in u_n C_\frac{\delta}{2}\right) =\PP(R>u_n)\PP\left(\thetav\in \Theta^\frac{\delta}{2}|R>u_n\right)\sim \PP(R>u_n)$
as $u\to\infty$,
\begin{eqnarray*} 
(1-\varepsilon')\sum_{n\geq 0} \PP\left(\X\in u_n C_\frac{\delta}{2}\right)
&\sim& (1-\varepsilon')\sum_{n\geq 0}\PP(R> u_n)\\
&\sim& \frac{1-\varepsilon'}{\|\cb\|+d\varepsilon}\int_u^\infty \PP(R>v)\ud v, 
\end{eqnarray*}
where $1-\varepsilon':=\frac{1-\delta'}{1+\delta''}$. The last asymptotic behaviour requires the subexponentiality of R and Assumption \ref{as0inttail}: for $n\leq n_1$ the term $dK$ in $u_n$ can be omitted by the long tail property and the assumption on the integrated tail distribution of $R$ ensures that the part of the sum is very large does not dominate the whole sum.
Finally, letting $\varepsilon,\varepsilon' \to 0$ yields the asymptotic lower bound of the ruin probability $\PP(\Sum_n \in uC \textrm{ for some } n\geq 1)$.

\subsubsection{Upper bound}
The upper bound 
\begin{displaymath}
\limsup_{u\to\infty} \frac{\PP(\Sum_n\in uC_\delta \textrm{ for some }n\geq 1)}{\frac{1}{\|\cb\|}\int_u^\infty \PP(R>v)\ud v}\leq 1
\end{displaymath}
follows from the deduction of Theorem 5.2 in \cite{Samorodnitsky} with small modifications. 
The relation $uC_\delta \subset uA$ yields 
\begin{eqnarray*} 
\PP(\Sum_n \in uC_\delta \textrm { for some } n\geq 1) &\leq& \PP(\Sum_n\in uA \textrm{ for some }n\geq 1) \\ 
&\sim& \int_0^\infty G_{\Z}(uA+v\cb)\ud v
\end{eqnarray*}
by Remark \ref{rem0samsun}, where $G_{\Z}$ denotes the distribution function of the positive random vector $\Z_1$ introduced in the remark. Assumption \ref{as0subexp}, the subexponentiality of the distribution $G_\X$ implies the multivariate long-tail property $G_\Z(uA+\ab)\sim G_\Z(uA)$ as $u\to\infty$ for any $\ab\in \RR^d$ of the distribution $G_{\Z}$ since
\begin{eqnarray*}
1\geq\lim_{u\to\infty} \frac{G_{\Z}(uA+v\cb +y\p)}{G_{\Z}(uA+v\cb)}
=\lim_{u\to\infty} \frac{G_{\X}(uA+v\cb)}{G_{\Z}(uA+v\cb)} \lim_{u\to\infty} \frac{G_{\Z}(uA+v\cb +y\p)}{G_{\X}(uA+v\cb +y\p)}=1.
\end{eqnarray*}
Next, we show the asymptotic equivalence $\int_0^\infty G_{\Z}(uA+v\cb)\ud v \sim \int_0^\infty G_{\X}(uA+v\cb)\ud v$ where $G_{\X}$ denotes the distribution function of the random vector $\X=\Z_1-Y_1\p$. Due to the fact that $A$ is an increasing set,
\begin{eqnarray*}
\int_0^\infty G_{\X}(uA+v\cb)\ud v  
\leq \int_0^\infty \PP(\Z_1\in uA+v\cb) \ud v =\int_0^\infty G_{\Z}(uA+v\cb)\ud v.
\end{eqnarray*}
Throughout the following calculations, $B$ denotes the distribution of the random variable $Y$.
By Fubini, Fatou lemma and the long-tail property of the distribution $H(u)$
\begin{eqnarray*} 
\lim_{u\to\infty} \int_0^\infty G_{\Z}(uA+v\cb)\ud v 
&=& \lim_{u\to\infty} \int_0^\infty \int_0^\infty G_{\X}(uA+v\cb- y\p)\ud B(y)\ud v\\
&\leq& \int_0^\infty \limsup_{u\to\infty}  \int_0^\infty G_{\X}(uA+v\cb- y\p)\ud v\ud B(y)\\
&=& \lim_{u\to\infty}  \int_0^\infty G_{\X}((u+v\|\cb\|)A)\ud v,
\end{eqnarray*}
which implies the asymptotic equivalence between the integrals as $u\to \infty$. Due to Assumption \ref{as0subexp}, also $\max\left(0,1-\int_0^\infty G_{\Z}(uA+v\cb)\ud v\right)$ is a subexponential distribution. Remark \ref{rem0samsun} together with the tail equivalence of the integrals implies
\begin{displaymath}
\lim_{u\to\infty} \frac{\PP\left(\sum_{i=1}^n\X_i\in uA \textrm{ for some }n\geq 1\right) }{\int_0^\infty G_{\X}(uA+v\cb) \ud v} =1
\end{displaymath}
and
\begin{eqnarray*} 
&&\int_0^\infty \PP(R>u+v\|\cb\|)\ud v \geq \int_0^\infty G_{\X}(uA+v\cb)\ud v \\
&\geq& \int_0^\infty \PP(R>u+v\|\cb\|)\PP\left(\thetav\in \Theta^\delta | R>u+v\|\cb\|\right)\ud v \\
&\geq& (1-\varepsilon)\int_0^\infty \PP(R>u+v\|\cb\|)\ud v \sim \frac{1}{\|\cb\|}\int_u^\infty \PP(R>v)\ud v
\end{eqnarray*}
shows the asymptotic upper bound.
\Endproof

\begin{rem}
In the proof of the upper bound Assumption \ref{as0inttail} is not needed. However, proving the lower bound using Lemma \ref{lem0geomapulem} requires the additional assumption that is fulfiled by a large class of distributions including Weibull distributions and lognormal distributions. Although the assumption does not hold for regularly varying random variables, a similar result that is based on a large deviations principle holds for them, see Theorem 3.1 in \cite{Hult1}. 
\end{rem}

\subsection{Proof of Theorem \ref{lem0equivalence}}
Due to Theorem \ref{lem0thetavakio}, it is enough to show that, under the additional assumption, Equation (\ref{eq0ruinprob}) implies Condition (\ref{eq0thetav}). Assuming the contrary of (\ref{eq0thetav}) we show a contradiction.
Choose $\Theta=\Theta_A$ in (\ref{eq0ruinprob}) and for all $0<\varepsilon\leq \delta$
\begin{eqnarray*}
\lim_{h\to\infty}\PP\left(\thetav \in \Theta_A^\varepsilon  |R>h\right) = p < 1 
\quad\textrm{and}\quad
\lim_{h\to\infty}\PP\left(\thetav \in \Theta_B^\varepsilon  |R>h\right) = 1-p <1,
\end{eqnarray*}
where $\delta>0$ is such that $\Theta_A^\delta\cap\Theta_B^\delta = \emptyset$.
We denote the set $uC_{\Theta_A,\delta}$ defined in (\ref{eq0cdelta}) by $uC_A$ and similarly $uC_B=uC_{\Theta_B,\delta}$.
The proof of Theorem \ref{lem0thetavakio} implies the lower bounds
\begin{equation} \label{lowbou1} 
\liminf_{u\to\infty} \frac{\PP(\Sum_n\in uC_{A} \textrm{ for some }n)}{\frac{1}{\|\cb\|}\int_u^\infty \PP(R>v) \ud v}\geq p 
\end{equation}
and
\begin{equation} \label{lowbou2}
\liminf_{u\to\infty} \frac{\PP(\Sum_n\in uC_{B} \textrm{ for some }n)}{\frac{1}{\|\cb\|}\int_u^\infty \PP(R>v) \ud v}\geq 1-p.
\end{equation}
We will show that the corresponding upper bounds are less than $1$. We write
\begin{eqnarray} \label{eq0probunion}
1&=&\lim_{u\to\infty} \frac{\PP(\Sum_n\in u(C_{A}\cup C_{B}) \textrm{ for some }n)}{\frac{1}{\|\cb\|}\int_u^\infty \PP(R>v) \ud v}\nonumber\\
&=&\lim_{u\to\infty} \frac{\PP(\Sum_n\in uC_{A} \textrm{ for some }n)+\PP(\Sum_n\in uC_{B} \textrm{ for some }n)}{\frac{1}{\|\cb\|}\int_u^\infty \PP(R>v) \ud v}\nonumber\\
&&\quad -\lim_{u\to\infty} \frac{\PP(\Sum_n\in uC_{A} \textrm{ for some }n,\Sum_k\in uC_{B} \textrm{ for some }k)}{\frac{1}{\|\cb\|}\int_u^\infty \PP(R>v) \ud v}.
\end{eqnarray}
The existence of the limits follows once we showed that the last term where the random walk visits both sets is zero. To do this, we sum over all possible times when the random walk visits the first set and the number of steps that it needs to reach the second set. Now,
\begin{eqnarray*} 
&&\PP(\Sum_n\in uC_{A} \textrm{ for some }n,\Sum_k\in uC_{B} \textrm{ for some }k)\\
&\leq&  \PP(\Sum_n\in uC_A, \Sum_{n+k}\in uC_B \textrm{ for some }n,k \geq 1) \\
&&+  \PP(\Sum_n\in uC_B, \Sum_{n+k}\in uC_A \textrm{ for some }n,k\geq 1).
\end{eqnarray*}
Below, we consider only the first probability where the random walk hits at first the set $uC_A$ and after $k$ steps the set $uC_B$. The calculations for the second probability are similar. 
We divide the probability into two parts depending on the number of steps $k$ that are needed to reach the second set. In order for the process to move from the first set to the second set in $k$ steps for $k\leq k(u)$ it has to hold $\|\Sum_k\| > u\ \dist(C_A,C_B)$, where $\dist(\cdot,\cdot)$ denotes the distance between two sets in $L_1$-norm and thus
\begin{eqnarray*} 
&&\PP(\Sum_{n+k}\in uC_B \textrm{ for some }k\leq k(u) \ |\ \Sum_{n}\in uC_A \textrm{ for some }n\geq 1 )\\
&\leq& \sum_{k=1}^{k(u)} \PP(\|\Sum_k\| > u\ \dist(C_A,C_B)) \\
&\leq& \sum_{k=1}^{k(u)}\PP\left(\sum_{i=1}^k R_i >u\ \dist(C_A,C_B)\right)
\sim (k(u))^2\PP(R>u \ \dist(C_A,C_B)).
\end{eqnarray*}
The asymptotic equivalence follows from the subexponentiality of the random variable $R$. Taking $k(u)=o\left(\sqrt{\left(\hanta(\dist(C_A,C_B)u)\right)^{-1}}\right)$ implies $(k(u))^2\PP(R>\dist(C_A,C_B)u)=o(1)$.

For $k>k(u)$, we distinguish between the cases whether the set $uC_B$ can be reached by drifting into the direction of the expectation or not. Therefore, we set $\gamma := \inf\{t:((t\Theta_A^\delta -\cb^{\varepsilon}) \cap C_B) \neq \emptyset \}>1$, where 
\begin{displaymath}
-\cb^{\varepsilon}:= \nolla \cup
\left\{-\x:\left\|\frac{\x}{\|\x\|}-\frac{\cb}{\|\cb\|}\right\|<\varepsilon\right\}
\end{displaymath}
is the smallest cone that contains $\cup_{k\geq 1} B(-k\cb,k\varepsilon)$. If $\gamma=\infty$, a drift from the set $uC_A$ into the set $uC_B$ is not possible so
\begin{eqnarray*}
&&  \PP(\Sum_n\in uC_A \textrm{ for some }n\geq 1, \Sum_{n+k} \in uC_B \textrm{ for some } k>k(u))\\
&\leq& \PP(\Sum_n \in uC_A \textrm{ for some }n\geq 1) \PP\left(\Sum_{k}\in u(C_B-C_A)\textrm{ for some }k>k(u)\right)
\end{eqnarray*}
and
\begin{eqnarray} \label{prob0erotus}
1\geq 
\sum_{k>k(u)} \PP(\Sum_k \in u(C_B-C_A), \Sum_j\notin u(C_B-C_A), j=1,\dots,k-1).
\end{eqnarray}
Since each of the terms in the sum is bounded by $\PP(\Sum_k\in u(C_B-C_A))\leq 1-\PP(\Sum_k\in \cb^{\varepsilon}) \underset{k \to\infty}{\longrightarrow} 0$ and $k(u)\to\infty$ as $u\to \infty$, the probability (\ref{prob0erotus}) converges to zero as $u$ grows.

If $\gamma<\infty$, drifting from $uC_A$ into $uC_B$ is possible but only if $\Sum_n\in \gamma uC_A$. In the case where $\Sum_n\in uC_A\backslash \gamma uC_A$, drifting is not possible and the above calculations hold, if we exchange the set $u(C_B-C_A)$ by $u(C_B-C_A\backslash \gamma C_A)$. 
In the case where $\Sum_n\in \gamma uC_A$, we divide the set $\gamma uC_A$ into smaller sets and condition on the fact that $\Sum_n$ belong to the smaller sets. Using the law of total probability,
\begin{eqnarray} \label{prob0erotus2} 
&&\PP(\Sum_n\in uC_A \textrm{ for some }n,\ \Sum_{n+k}\in uC_B \textrm{ for some }k>k(u))\\
&=&\sum_{l=1}^\infty \sum_{n=1}^\infty \PP(\Sum_n \in l\gamma uC_A\backslash (l+1)\gamma uC_A,\  \Sum_j \notin uC_A, j=1\dots,n-1)\times\ldots\nonumber\\
&&\qquad \ldots\times\PP(\Sum_{n+k} \in uC_B \textrm{ for some }k>k(u) \ |\ \Sum_n\in l\gamma uC_A\backslash (l+1)\gamma uC_A)\nonumber\\
&\leq& \sum_{l=1}^\infty \PP(\Sum_n\in l\gamma uC_A\backslash (l+1)\gamma uC_A \textrm{ for some }n)\times \ldots\nonumber\\
&& \qquad \ldots\times\PP(\|\Sum_k\|\leq c_l u\textrm{ for some }k>k(u))\label{ineq03}\\
&=& \sum_{k>k(u)}\sum_{l=1}^\infty \PP(\Sum_n\in l\gamma uC_A\backslash (l+1)\gamma uC_A \textrm{ for some }n)\times \ldots\nonumber\\
&& \qquad \ldots\times\PP(\|\Sum_k\|\leq c_l u, \|\Sum_j\|> c_l u, j=k(u)+1,\dots, k-1)\nonumber\\
&\leq& \PP(\Sum_n\in uC_A\textrm{ for some }n)\ \leq\  1 \nonumber
\end{eqnarray}
where 
\begin{displaymath}
c_l:= \sup\left\{t>0:\{(\gamma lC_A\backslash \gamma (l+1)C_A -t\cb^\varepsilon) \cap C_B\}\neq \emptyset\right\} 
\end{displaymath}
with $-t\cb^\varepsilon:=\{\x\in-\cb^\varepsilon:\|\x\|>t\}$ describes the largest distance between the sets $\gamma lC_A\backslash \gamma (l+1)C_A$ and $C_B$ in the direction of the vector $\cb$. 
Now, for every fixed $l$
\begin{eqnarray*} 
&&\PP(\|\Sum_k\| \leq c_lu, \|\Sum_j\|> c_lu, j=k(u)+1,\dots,k-1)\\
&\leq& \PP(\|\Sum_k\| \leq c_lu)
\leq 1-\PP(\Sum_k\in B(-k\cb,k\varepsilon)) \underset{u\to\infty}{\longrightarrow} 0
\end{eqnarray*}
due to the law of large numbers. The last inequality follows from the assumption $u^2\hanta(u)=o(1)$, since one can choose the function $k(u)$ such that $u=o(k(u))$. Due to the fact that the sum (\ref{ineq03}) is finite, Probability (\ref{prob0erotus2}) converges to zero as $u$ grows.

Similar calculations for the case where the random walk visits at first the set $uC_B$ and then the set $uC_A$ finally yield that the last term on the right-hand side in (\ref{eq0probunion}) is zero which implies that the lower bounds (\ref{lowbou1}) and (\ref{lowbou2}) are also upper bounds. Thus, 
\begin{displaymath}
\limsup_{u\to\infty}\PP(\Sum_n \in uC_{A} \textrm{ for some }n\geq 1) <1
\end{displaymath}
which yields the contradiction.
\Endproof

\section{Applications and Examples}
The class of distributions that fulfil Assumptions \ref{as0r} - \ref{as0subexp} contains distributions that are subexponential, but have a lighter tail than the class of multivariate regularly varying distributions. In particular, the components of such distributions can have finite moments of all orders.

\begin{rem} \label{rem0inttail}
The most common subexponential distributions such as regularly varying distributions, Weibull distributions with $\beta\in (0,1)$, the Lognormal distribution, Benktander distributions of type  I and II as well as the Loggamma distributions have the property that also their integrated tail is subexponentially distributed, see \cite{Embrechts1}. 
\end{rem}

\begin{exa}
Suppose $\X=R\thetav$ where $R$ has a common subexponential distribution as in Remark \ref{rem0inttail} and Assumptions \ref{as0thetav} and \ref{as0expectation} hold. Assumption \ref{as0thetav} implies
\begin{eqnarray*}
\PP(\X \in uA+v\cb) 
&=& \PP(R>u+v\|\cb\|) \PP\left(\thetav \in \RR^d_+|R>u+v\|\cb\|\right)\\
&\sim& \PP(R>u+v\|\cb\|).
\end{eqnarray*}
Thus,
\begin{equation} \label{eq0aequcdelta}
\int_0^\infty \PP(R\thetav \in uA+v\cb) \ud v \sim
\frac{1}{\|\cb\|}\int_u^\infty \PP(R>v) \ud v
\end{equation}
is finite and the integral on the right-hand side is subexponential. Due to the tail equivalence of subexponential distributions, the distribution $H(u)$ is subexponential as well so Assumptions \ref{as0finite} and \ref{as0subexp} hold. 
\end{exa}

\begin{exa}
If $R$ has lognormal distribution or Weibull distribution with parameter $0<\beta<1$, Assumption \ref{as0inttail} holds always. Furthermore, these distribution fulfil the assumption $u^2\hanta(u)=o(1)$ of Theorem \ref{lem0equivalence}. 
If $R$ is Weibull distributed with parameter $0<\beta<1$, a suitable function $k(u)$ is $k(u)= \exp(u^\beta/2-\varepsilon)$ for some $\varepsilon>0$. If $R$ follows a lognormal distribution, $k(u)=u^2$ can be used. 
\end{exa}

\begin{rem}
Under Assumptions \ref{as0r} - \ref{as0subexp},
\begin{displaymath}
\PP(\Sum_n\in uA \textrm{ for some } n\geq 1) \sim \PP(\Sum_n \in uC_\delta \textrm{ for some }n\geq 1) \qquad \textrm{ as } u\to\infty
\end{displaymath}
due to the asymptotic relation (\ref{eq0aequcdelta}). Interpreting $\Sum_n$ as the net-payout process of an insurance company with different lines of business, the total net-payout of the company $\sum_{k=1}^d S_n^k$ is large because ruin occurs in the direction of the set $\Theta^\delta$.
\end{rem}

\begin{rem}
Extending the bivariate model to different subexponential distributions in different directions, one can recognize some hidden subexponentiality analogue to hidden regular variation. Similarly as in \cite{Das1}, hidden subexponentiality is applicable whenever $\Theta$ is a union of sets with different distributions for $R$. For instance, in the case where different sets have Weibull distributed vector lengths with different parameters, the ruin probability is determined by the distribution with the smallest parameter and the other parameters are reflected in the hidden part.
\end{rem}

\section*{Acknowledgements}
The financial support from the Doctoral Programme of Mathematics and Statistics of the University of Helsinki, Domast, is gratefully acknowledged.
Special thanks are due to Jaakko Lehtomaa for supportive discussions during the writing of this paper. Suggestions made by an anonymous referee greatly improved the manuscript. Declaration of interest: none.

\bibliography{references3}
\bibliographystyle{abbrv}

\end{document}